\documentclass[pamm,a4paper,fleqn]{w-art}
\usepackage{times,cite,w-thm}
\usepackage[T1]{fontenc}
\usepackage[utf8]{inputenc}
\usepackage{amsmath,bm,mathtools,amssymb}
\usepackage{graphicx}
\usepackage{booktabs}
\usepackage{multirow}
\usepackage{relsize}

\usepackage{tikz}
\usepackage{pgfplots, pgfplotstable}
\pgfplotsset{compat=newest}
\usetikzlibrary{external}
\usepackage{siunitx}
\newcommand{\mR}{\mathbb{R}}

\newcommand{\SO}{{SO}}
\let\so\undefined
\newcommand{\so}{\mathfrak{so}}
\newcommand{\SE}{{S\!E}}

\newcommand{\Skw}{\operatorname{Skw}}

\newcommand{\T}{^{\!\mathrm{T}}}

\newcommand{\diff}[1][]{\mathrm{d}#1}

\newcommand{\vga}{\bm{\gamma}}

\newcommand{\vth}{\bm{\theta}}

\newcommand{\vka}{\bm{\kappa}}

\newcommand{\vph}{\bm{\phi}} 

\newcommand{\vom}{\bm{\omega}}

\newcommand{\vb}{\mathbf b}
\newcommand{\vc}{\mathbf c}

\newcommand{\ve}{\mathbf e}
\newcommand{\vf}{\mathbf f}
\newcommand{\vg}{\mathbf g}

\newcommand{\vm}{\mathbf m}
\newcommand{\vn}{\mathbf n}

\newcommand{\vp}{\mathbf p}
\newcommand{\vq}{\mathbf q}
\newcommand{\vr}{\mathbf r}
\newcommand{\vs}{\mathbf s}

\newcommand{\vu}{\mathbf u}
\newcommand{\vv}{\mathbf v}

\newcommand{\vA}{\mathbf A}
\newcommand{\vB}{\mathbf B}
\newcommand{\vC}{\mathbf C}
\newcommand{\vD}{\mathbf D}

\newcommand{\vF}{\mathbf F}

\newcommand{\vH}{\mathbf H}
\newcommand{\vI}{\mathbf I}

\newcommand{\vK}{\mathbf K}

\newcommand{\vM}{\mathbf M}

\newcommand{\vP}{\mathbf P}
\newcommand{\vQ}{\mathbf Q}

\begin{document}

\TitleLanguage[EN]
\title[The short title]{Non-unit quaternion parametrization of a Petrov--Galerkin Cosserat rod finite element}

\author{\firstname{Jonas} \lastname{Harsch}\inst{1,}%
\footnote{Corresponding author: \ElectronicMail{harsch@inm.uni-stuttgart.de}}}
\address[\inst{1}]{\CountryCode[DE] Institute for Nonlinear Mechanics, University of Stuttgart, Stuttgart, Germany}
\author{\firstname{Simon R.} \lastname{Eugster}\inst{1}}

\AbstractLanguage[EN]
\begin{abstract}
The application of the Petrov--Galerkin projection method in Cosserat rod finite element formulations offers significant advantages in simplifying the expressions within the discrete virtual work functionals. Moreover, it enables a straight-forward and systematic exchange of the ansatz functions, specifically for centerline positions and cross-section orientations. In this concise communication, we present a total Lagrangian finite element formulation for Cosserat rods that attempts to come up with the least required concepts. The chosen discretization preserves objectivity and allows for large displacements/ rotations and for large strains. The orientation parametrization with non-unit quaternions results in a singularity-free formulation.
\end{abstract}
\maketitle                   

\section{Introduction}
This article complements the two papers \cite{Harsch2023a, Eugster2023a} on Petrov--Galerkin rod finite formulations for Cosserat rods. The cross-section orientations are parameterized using non-unit quaternions instead of total rotation vectors, which require additionally the concept of the complement rotation vector for a singularity-free parametrization. To keep the formulation as simple as possible, we opt for the $\mR^{12}$-interpolation for the ansatz functions, see \cite{Betsch2002, Romero2002, Eugster2023a}.

The paper is structured as follows. In Section~\ref{sec:Cosserat_rod_theory}, the Cosserat rod theory is recapitulated very briefly; mainly to introduce all quantities required for the further finite element formulation. For those interested in additional comments as well as a thorough introduction and explanation of the chosen notation, we recommend reading \cite{Harsch2023a, Eugster2023a}. The Petrov--Galerkin finite element formulation in terms of nodal non-unit quaternions is presented in Section~\ref{sec:PG_fem}. The last section on numerical experiments, investigates the static analysis of a helical spring in line with~\cite{Marino2017}. Additionally, the Wilberforce example from~\cite{Harsch2021a} with a helical spring with three coils is discussed.

\section{Cosserat rod theory}\label{sec:Cosserat_rod_theory}
Let $\xi \in \mathcal{J} = [0, 1] \subset \mR$ be the centerline parameter and let $t$ denote time. The motion of a Cosserat rod is captured by a time-dependent centerline curve represented in an inertial $I$-basis ${}_I\vr_{OP}={}_I\vr_{OP}(\xi, t) \in \mR^3$ augmented by the cross-section orientations $\vA_{IK}=\vA_{IK}(\xi, t) \in SO(3)=\{\vA \in \mR^{3 \times 3}| \vA\T \vA = \mathbf{1}_{3 \times 3} \wedge \det(\vA)= 1\}$. The subscripts $O$ and $P$ in the centerline curve refer to the origin and the centerline point, respectively. The cross-section orientation $\vA_{IK}$ can also be interpreted as a transformation matrix that relates the representation of a vector in the cross-section-fixed $K$-basis to its representation in the inertial $I$-basis.

The derivatives with respect to time $t$ and centerline parameter $\xi$ are denoted by $\dot{(\bullet)}$ and $(\bullet)_{,\xi}$, respectively. The variation of a function is indicated by $\delta(\bullet)$. With this, we can introduce the centerline velocity ${}_I \vv_{P}= \left({}_I \vr_{OP}\right)^\mathlarger\cdot$ and the virtual displacement ${}_I \delta \vr_{P} = \delta \left({}_I \vr_{OP}\right)$.
The angular velocity of the cross-section-fixed $K$-basis relative to the inertial $I$-basis, in components with respect to the $K$-basis, is defined by ${}_K \vom_{IK} \coloneqq j^{-1} \big(\vA_{IK}\T \left(\vA_{IK}\right)^\mathlarger\cdot \big)$,
where $j \colon \mR^3 \to \so(3) = \{\vB \in \mR^{3\times3} | \vB\T = -\vB\}$ is the linear and bijective map such that $\widetilde{\vom} \vr = j(\vom) \vr = \vom \times \vr$ for all $\vom, \vr \in \mR^3$. Analogously, the virtual rotations and the scaled curvature are defined as
${}_K \delta\vph_{IK} \coloneqq j^{-1} \big( \vA_{IK}\T \delta \left(\vA_{IK}\right) \big)$ and ${}_K \bar{\vka}_{IK} \coloneqq j^{-1} \big( \vA_{IK}\T (\vA_{IK})_{,\xi} \big)$, respectively. 
For the reference centerline curve ${}_I \vr_{OP}^0$, the length of the rod's tangent vector is $J = \|{}_I \vr_{OP, \xi}^0\|$. Thus, for a given centerline parameter $\xi$, the reference arc length increment is $\diff s = J \diff \xi$. The derivative with respect to the reference arc length $s$ of a function $\vf = \vf(\xi,t) \in \mR^3$ can then be defined as $\vf_{,s}(\xi,t) \coloneqq \vf_{,\xi}(\xi,t) /J(\xi)$.
The objective strain measures of a Cosserat rod are the curvature ${}_K \vka_{IK} = {}_K \bar{\vka}_{IK} / J$, which measures torsion and bending, together with the measures for dilatation and shear strains contained in ${}_K \vga = {}_K \bar{\vga} / J$ determined by ${}_K \bar{\vga} \coloneqq (\vA_{IK})\T {}_I \vr_{OP, \xi}$.

The internal virtual work of a Cosserat rod is defined as
\begin{equation}\label{eq:internal_virtual_work2}
\delta W^\mathrm{int} \coloneqq -\int_{\mathcal{J}} \big\{ 
({}_I \delta \vr_{P,\xi})\T \vA_{IK} {}_K \vn + ({}_K \delta \vph_{IK,\xi})\T  {}_K \vm 
- ({}_K \delta\vph_{IK})\T \left[ {}_K \bar{\vga} \times {}_K \vn + {}_K \bar{\vka}_{IK} \times {}_K \vm \right] 
\big\} \diff[\xi] \, ,
\end{equation}
where ${}_K \vn$ and ${}_K \vm$ denote the resultant contact forces and moments, respectively. For hyperelastic material models with a strain energy density with respect to the reference arc length $W = W({}_K \vga, {}_K \vka_{IK}; \xi)$, they can be determined by the constitutive relations ${}_K \vn = (\partial W / \partial {}_K \vga)\T$ and ${}_K \vm = (\partial W / \partial {}_K \vka_{IK})\T$.

Assume the line distributed external forces ${}_I \vb = {}_I \vb(\xi,t) \in \mR^3$ and moments ${}_K \vc ={}_K \vc(\xi,t) \in \mR^3$ to be given as densities with respect to the reference arc length. Moreover, for $i\in\{0,1\}$, point forces ${}_I \vb_i = {}_I \vb_i(t) \in \mR^3$ and point moments ${}_K \vc_i = {}_K \vc_i(t) \in \mR^3$ can be applied to the rod's boundaries at $\xi_0=0$ and $\xi_1=1$. The corresponding external virtual work functional is defined as
\begin{equation}\label{eq:external_virtual_work}
\delta W^\mathrm{ext} \coloneqq \int_{\mathcal{J}} \left\{ ({}_I\delta\vr_{P})\T {}_I \vb + ({}_K \delta \vph_{IK})\T {}_K \vc \right\} J \diff[\xi]
+ \sum_{i = 0}^1 \left[ ({}_I\delta\vr_{P})\T {}_I \vb_i + ({}_K \delta \vph_{IK})\T {}_K \vc_i \right]_{\xi_i} \, .
\end{equation}

In case ${}_I\vr_{OP}$ is the line of centroids, the inertial virtual work functional of the Cosserat rod can be written as
\begin{equation}\label{eq:inertia_virtual_work}
\delta W^\mathrm{dyn} \coloneqq -\int_{\mathcal{J}}
\big\{({}_I \delta \vr_{P} )\T A_{\rho_0} ({}_I\vv_p)^\mathlarger{\cdot} + ({}_K \delta \vph_{IK})\T
({}_K \vI_{\rho_0} ({}_K \vom_{IK})^\mathlarger{\cdot} + {}_K {\vom}_{IK} \times {}_K \vI_{\rho_0} {}_K \vom_{IK})\big\} J \diff[\xi]\, ,
\end{equation}
where $A_{\rho_0}$ is the cross-section mass density and ${}_K \vI_{\rho_0}$ the constant cross-section inertia tensor represented in the cross-section-fixed $K$-basis.

\section{Petrov--Galerkin finite element formulation}\label{sec:PG_fem}
The rod's parameter space $\mathcal{J}$ is divided into $n_\mathrm{el}$ linearly spaced element intervals $\mathcal{J}^e = [\xi^{e}, \xi^{e+1})$ via $\mathcal{J} = \bigcup_{e=0}^{n_\mathrm{el}-1} \mathcal{J}^e$. For a $p$-th order finite element, the closure of each of the intervals $\mathcal{J}^e$ contains $p + 1$ evenly spaced points $\xi^e_i \in \mathrm{cl}(\mathcal{J}^e) = [\xi^{e}, \xi^{e+1}]$ with $i \in \{0, \dots, p\}$ such that $\xi^e_0 = \xi^e < \xi^e_1 < \dots < \xi^e_p = \xi^{e+1}$. Note, for $e \in \{0, \ldots, n_{\mathrm{el}} -2 \}$, the points $\xi^e_p=\xi^{e+1}_0$ denote the same point $\xi^{e+1}$, which is the boundary point of the adjacent element intervals. It is convenient to use both indexations in the following. For a given element interval $\mathcal{J}^e = [\xi^e, \xi^{e+1})$, the $p$-th order Lagrange basis function and derivative of node $i\in \{0,\dots,p\}$ are
\begin{equation}\label{eq:Lagrangian_polynomials}
N^{p,e}_i(\xi) = \underset{\substack{0 \leq j \leq p \\ j\neq i}}{\prod} \frac{\xi - \xi^e_j}{ \xi^e_i - \xi^e_j} \quad \mathrm{and}
\quad N^{p,e}_{i,\xi}(\xi) = N_i^{p,e}(\xi) \underset{\substack{k=0 \\ k \neq i}}{\sum^{p}} \frac{1}{\xi - \xi^e_k} \, ,
\end{equation}
where $\xi^e_i$, $\xi^e_j$, and $\xi^e_k$ are the points contained in the set $\{\xi^e_0 = \xi^e, \xi^e_1, \dots, \xi^e_p = \xi^{e+1}\}$.

The centerline curve ${}_I \vr_{OP}$ and the cross-section orientations $\vA_{IK}$ are approximated by interpolating nodal centerline points ${}_I\vr_{OP^e_i}(t)\in \mR^3$ and nodal transformation matrices $\vA_{IK^e_i}(t)\in\SO(3)$. For each node $ i \in \{0,\dots,p\}$ within element $e \in \{0,\dots, n_\mathrm{el}-1\}$, it will hold that ${}_I\vr_{OP^e_i}(t) = {}_I \vr_{OP}(\xi^e_i,t)$ and $\vA_{IK^e_i}(t) = \vA_{IK}(\xi^e_i,t)$. In contrast to \cite{Harsch2023a, Eugster2023a}, the nodal transformation matrices 
\begin{equation}\label{eq:quat2rot}
\vA_{IK^e_i} = \vA(\vP^e_i) = \mathbf{1}_{3 \times 3} + 2 \left((\widetilde{\vp}^e_i)^2 + p^e_{0,i} \, \widetilde{\vp}^e_i\right) / \|\vP_i^e\|^2
\end{equation}
are parametrized by nodal non-unit quaternions  $\vP^e_i(t) = (p^e_{0,i}(t), \vp^e_i(t)) \in \mR^4$ with the scalar part $p^e_{0,i}(t) \in \mR$ and the vectorial part $\vp^e_i(t) \in \mR^3$, see \cite{Rucker2018}. Note that \eqref{eq:quat2rot} is formulated in such a way to return orthogonal matrices also for non-unit quaternions.

Accordingly, the $N=(p n_\mathrm{el} + 1)$ nodal generalized position coordinates $\vq^e_i(t) = ({}_I \vr_{OP^e_i}, \vP^e_i)(t) \in \mR^7$ are given by the nodal centerline points ${}_I \vr_{OP^e_i}$ and the nodal non-unit quaternions $\vP^e_i$ resulting in $n_{\vq} = 7N$ positional degrees of freedom for the discretized rod. The nodal quantities can be assembled in the global tuple of generalized position coordinates $\vq(t) = \big(\vq^0_0, \dots, \vq^0_{p-1}, \dots, \vq^e_{0}, \dots, \vq^e_{p-1}, \dots, \vq^{n_\mathrm{el}-1}_0, \dots,\vq^{n_\mathrm{el}-1}_{p-1}, \vq^{n_\mathrm{el}-1}_p\big)(t) \in \mR^{n_{\vq}}$. For $e \in \{0, \ldots, n_{\mathrm{el}} -2 \}$, the coordinates $\vq^e_p=\vq^{e+1}_0$ refer to the same nodal coordinates. Introducing an appropriate Boolean connectivity matrix $\vC_e \in \mR^{7(p+1) \times n_{\vq}}$, the element generalized position coordinates $\vq^e(t) = \big(\vq^e_0, \dots, \vq^e_p\big)(t) \in \mR^{7(p+1)}$ can be extracted from $\vq$ via $\vq^e = \vC_{e} \vq$. Note that during a numerical implementation it is advisable to slice arrays instead of multiply them with Boolean matrices.

In the sense of \cite{Betsch2002, Romero2002}, both the nodal centerline points and the cross-section orientations are interpolated by $p$-th order Lagrangian polynomials. Using the characteristic function $\chi_{\mathcal{J}^e} \colon \mathcal{J} \to \{0, 1\}$, which is one for $\xi \in \mathcal{J}^e = [\xi^e, \xi^{e+1})$ and zero elsewhere, together with the $p$-th order Lagrange basis functions \eqref{eq:Lagrangian_polynomials}, the ansatz functions for centerline and cross-section orientations are
\begin{equation}\label{eq:R12_interpolation}
{}_I \vr_{OP}(\xi, \vq) = \sum_{e=0}^{n_\mathrm{el}-1} \chi_{\mathcal{J}^e}(\xi) \sum_{i=0}^p 
N^{p,e}_i(\xi) {}_I \vr_{OP^e_i}
\quad \mathrm{and} \quad
\vA_{IK}(\xi, \vq) = \sum_{e=0}^{n_\mathrm{el}-1} \chi_{\mathcal{J}^e}(\xi)
\sum_{i=0}^p 
N^{p,e}_i(\xi) \vA(\vP^e_i)\, .
\end{equation}

The discretized version of the curvature strain is computed as
\begin{equation}\label{eq:discrete_curvature}
{}_K \vka_{IK} = j^{-1}\big(\Skw(\vA_{IK}\T \vA_{IK,\xi})\big) / J \, ,
\end{equation}
where the map $\Skw(\vM) = \frac{1}{2}(\vM - \vM\T) \in \so(3)$ extracts the skew-symmetric part of the matrix $\vM\in\mR^{3\times3}$. Hence, the curvature can efficiently be computed using $j^{-1}(\Skw(\vM)) = \tfrac{1}{2}(M_{32} - M_{23}, \ M_{13} - M_{31}, \ M_{21} - M_{12})$. 

At the same $N$ nodes as for the nodal generalized position coordinates, we introduce the nodal generalized virtual displacements $\delta \vs^e_i(t) = ({}_I \delta \vr_{P^e_i}, {}_{K^e_i} \delta \vph_{IK^e_i})(t) \in \mR^6$ given by the nodal virtual centerline displacement ${}_I \delta \vr_{P^e_i}(t) \in \mR^3$ and the nodal virtual rotation ${}_{K^e_i} \delta \vph_{IK^e_i}(t) \in \mR^3$. In analogy to the generalized virtual displacements, we also introduce the nodal generalized velocities $\vu^e_i(t) = ({}_I \vv_{P^e_i}, {}_{K^e_i} \vom_{IK^e_i})(t) \in \mR^6$ given by the nodal centerline velocity ${}_I \vv_{P^e_i}(t) \in \mR^3$ and the nodal angular velocity ${}_{K^e_i} \vom_{IK^e_i}(t) \in \mR^3$. Similar to the generalized position coordinates $\vq$, the nodal generalized virtual displacements and velocities are assembled in the global tuple of generalized virtual displacements $\delta \vs(t) \in \mR^{n_{\vu}}$ and velocities $\vu(t) \in \mR^{n_{\vu}}$. In contrast to the nodal position coordinates, there are only six nodal generalized virtual displacements or velocity coordinates resulting in $n_{\vu} = 6N$ generalized virtual displacements or velocity degrees of freedom for the discretized rod. 

Consequently, we require a new Boolean connectivity matrix $\vC_{\vu, e}\in \mR^{6(p+1) \times n_{\vu}}$, which extracts the element generalized virtual displacements  $\delta \vs^e(t) = (\delta \vs^e_0, \dots, \delta \vs^e_p)(t) \in \mR^{6(p+1)}$ and  velocities $\vu^e(t) = (\vu^e_0, \dots, \vu^e_p)(t) \in \mR^{6(p+1)}$ from the global quantities via $\delta \vs^e = \vC_{\vu,e} \delta \vs$ and $\vu^e = \vC_{\vu,e} \vu$. By further introducing the Boolean connectivity matrices $\vC_{\vr, i} \in \mR^{3 \times 6(p+1)}$, the nodal virtual centerline displacements ${}_I \delta \vr_{P^e_i}$ and centerline velocities ${}_I \vv_{P^e_i}$ can be extracted from the element generalized virtual displacements $\delta \vs^e$ and velocities $\vu^e$ via ${}_I \delta \vr_{P^e_i} = \vC_{\vr, i} \delta \vs^e$ and ${}_I \vv_{P^e_i} = \vC_{\vr, i} \vu^e$, respectively. Identical extraction operations hold for the nodal virtual rotations ${}_{K^e_i} \delta \vph_{IK^e_i} = \vC_{\vph,i} \delta \vs^e$ and angular velocities ${}_{K^e_i} \vom_{IK^e_i}= \vC_{\vph,i} \vu^e$, where $\vC_{\vph, i} \in \mR^{3 \times 6(p+1)}$. The test functions are then given by interpolating the nodal generalized virtual displacements by $p$-th order Lagrangian basis functions \eqref{eq:Lagrangian_polynomials} in agreement with
\begin{equation}\label{eq:virtual_displacement_interpol}
{}_I \delta \vr_{P}(\xi, \delta \vs) = \sum_{e=0}^{n_\mathrm{el} - 1} \chi_{\mathcal{J}^e}(\xi) \sum_{i=0}^p N^{p,e}_i(\xi) {}_I \delta \vr_{P^e_i} \quad \mathrm{and} \quad
{}_K \delta \vph_{IK}(\xi, \delta \vs) = \sum_{e=0}^{n_\mathrm{el} - 1} \chi_{\mathcal{J}^e}(\xi)\sum_{i=0}^p N^{p,e}_i(\xi) {}_{K^e_i} \delta \vph_{IK^e_i} \, .
\end{equation}
Note that the interpolation of the virtual rotations must be understood in the sense of a Petrov--Galerkin projection, where the virtual rotations are not obtained from a consistent variation of the ansatz functions~\eqref{eq:R12_interpolation}.

To obtain a constant and symmetric mass matrix in the discretized formulation, see~\eqref{eq:discretized_inertial_virtual_work} below, the velocities are considered as independent fields and are interpolated with the same interpolation as the virtual displacements and rotations as
\begin{equation}\label{eq:velocity_interpolation}
{}_I \vv_{P}(\xi, \vu) = \sum_{e=0}^{n_\mathrm{el} - 1} \chi_{\mathcal{J}^e}(\xi)\sum_{i=0}^p N^{p,e}_i(\xi) {}_I \vv_{P^e_i} \quad \mathrm{and} \quad
{}_K \vom_{IK}(\xi, \vu) = \sum_{e=0}^{n_\mathrm{el} - 1} \chi_{\mathcal{J}^e}(\xi) \sum_{i=0}^p N^{p,e}_i(\xi) {}_{K^e_i} \vom_{IK^e_i} \, .
\end{equation}
The independent introduction of velocity fields \eqref{eq:velocity_interpolation} demands an additional relation defining the coupling between position coordinates $\vq$ and velocity coordinates $\vu$. This coupling is given by the nodal kinematic differential equations
\begin{equation}\label{eq:nodal_kinematic_equation}
\dot{\vq}^e_i =
\begin{pmatrix}
{}_I\dot{\vr}_{OP^e_i} \\
\dot{\vP}^e_i
\end{pmatrix} =
\begin{pmatrix}
\mathbf{1}_{3 \times 3} & \mathbf{0}_{3 \times 3} \\
\mathbf{0}_{4 \times 3} & \vQ(\vP^e_i)
\end{pmatrix}
\begin{pmatrix}
{}_I \vv_{P^e_i} \\
{}_{K^e_i} \vom_{IK^e_i}
\end{pmatrix} =
\vF(\vq^e_i) \vu^e_i \, , \quad \textrm{where} \; \,
\vQ(\vP) = \frac{1}{2}
\begin{pmatrix}
-\vp\T \\
p_0 \mathbf{1}_{3 \times 3} + \widetilde{\vp}
\end{pmatrix} \, ,
\end{equation}
cf.~\cite{Rucker2018}. The nodal kinematic equations~\eqref{eq:nodal_kinematic_equation} can easily be assembled to a global kinematic differential equation of the form $\dot{\vq} = \vB(\vq) \vu$. Note that the kinematic differential equation is linear in $\vq$ too. This allows to write the relation also in the form $\dot{\vq} = \vD(\vu) \vq$, see \cite{Rucker2018} for more details.

Inserting the test functions~\eqref{eq:virtual_displacement_interpol} together with the corresponding approximations for centerline, cross-section orientations~\eqref{eq:R12_interpolation} and strain measures into \eqref{eq:internal_virtual_work2}, the continuous internal virtual work is approximated by $\delta W^\mathrm{int}(\vq; \delta \vs) = \delta \vs\T \vf^{\mathrm{int}}(\vq)$, where the internal generalized forces are computed element-wise by
\begin{equation}\label{eq:internal_gen_forces}
\begin{aligned}
\vf^{\mathrm{int}}(\vq) & = \sum_{e=0}^{n_\mathrm{el} - 1} \vC_{\vu, e}\T \vf^{\mathrm{int}}_e(\vC_e \vq) \, , \\
\vf^{\mathrm{int}}_e(\vq^e) & = -\int_{\mathcal{J}^e} \sum_{i=0}^{p}\Big\{ N^{p,e}_{i,\xi} \vC_{\vr, i}\T \vA_{IK} {}_K \vn + N^{p,e}_{i,\xi} \vC_{\vph, i}\T {}_K \vm 
-N^{p,e}_{i} \vC_{\vph, i}\T \left({}_K \bar{\vga}\times {}_K \vn + {}_K\bar{\vka}_{IK} \times {}_K \vm \right) \Big\} \diff[\xi] \, .
\end{aligned}
\end{equation}
Similarly, the external virtual work~\eqref{eq:external_virtual_work} is discretized by $\delta W^\mathrm{ext}(t, \vq; \delta \vs) = \delta \vs\T \vf^{\mathrm{ext}}(t, \vq)$ with
\begin{align}
\vf^{\mathrm{ext}}(t, \vq) & = \sum_{e=0}^{n_\mathrm{el} - 1} \vC_{\vu,e}\T \vf^{\mathrm{ext}}_e(t, \vC_e \vq) + \vC_{\vu, 0}\T \left[\vC_{\vr, 0}\T {}_I \vb_0 {+} \vC_{\vph, 0}\T {}_K \vc_0 \right]_{\xi=0} +
\vC_{\vu, n_\mathrm{el} - 1}\T \left[\vC_{\vr,p}\T {}_I \vb_1 {+} \vC_{\vph, p}\T {}_K \vc_1 \right]_{\xi=1}\, , \nonumber \\
\vf^{\mathrm{ext}}_e(t, \vq^e) & = \int_{\mathcal{J}^e} \sum_{i=0}^{p}\Big\{ N^{p,e}_{i}  \vC_{\vr, i}\T {}_I \vb + N^{p,e}_{i} \vC_{\vph, i}\T {}_K \vc \Big\} J \diff[\xi]
 \, .
\end{align}
Finally, inserting  \eqref{eq:virtual_displacement_interpol} and \eqref{eq:velocity_interpolation} into the inertial virtual work functional~\eqref{eq:inertia_virtual_work} yields the discrete counterpart $\delta W^\mathrm{dyn}(\vu;\delta \vs) = -\delta \vs\T \big(\vM \dot{\vu} + \vf^{\mathrm{gyr}}(\vu) \big)$, 
where we have introduced the symmetric and constant mass matrix
\begin{equation}\label{eq:discretized_inertial_virtual_work}
\vM= \sum_{e=0}^{n_\mathrm{el} - 1} \vC_{\vu,e}\T \vM_e \vC_{\vu,e} \, , \quad
\vM_e = 
\int_{\mathcal{J}^e} \sum_{i=0}^{p} \sum_{k=0}^{p} N^{p,e}_{i} N^{p,e}_{k} \Big\{
A_{\rho_0} \vC_{\vr, i}\T  \vC_{\vr, k} 
+ \vC_{\vph, i}\T {}_K \vI_{\rho_0} \vC_{\vph, k}
\Big\} J \diff[\xi] \, ,
\end{equation}
and the gyroscopic forces
\begin{equation}
\vf^{\mathrm{gyr}}(\vu) = \sum_{e=0}^{n_\mathrm{el}-1} \vC_{\vu,e}\T \vf_e^{\mathrm{gyr}}(\vC_{\vu,e} \vu) \, , \quad
\vf^{\mathrm{gyr}}_e(\vu^e) = \int_{\mathcal{J}^e} \sum_{i=0}^{p} N^{p,e}_i \Big\{\vC_{\vph, i}\T ({}_K {\vom}_{IK} \times {}_K \vI_{\rho_0} {}_K \vom_{IK}) \Big\} J \diff[\xi]  \, .
\end{equation}
Element integrals of the form $\int_{\mathcal{J}^e} f(\xi) \diff[\xi]$ arising in the discretized external and gyroscopic forces, as well as in the mass matrix, are subsequently computed using a Gauss--Legendre quadrature rule with $\operatorname{ceil}[(p + 1)^2 / 2]$ quadrature points. To alleviate locking, the internal generalized forces \eqref{eq:internal_gen_forces} are integrated by a reduced $p$-point quadrature rule.

Applying the principle of virtual work, which requires the total virtual work functional to vanish, we readily obtain the system dynamics in the form 
\begin{equation}
\begin{aligned}
\dot{\vq} &= \vB(\vq) \vu \, , \\
\dot{\vu} &= \vM^{-1} \left(\vf^\mathrm{gyr}(\vu) + \vf^\mathrm{int}(\vq) + \vf^\mathrm{ext}(t, \vq)\right) \, ,
\end{aligned}
\end{equation}
where the two lines correspond to the global kinematic differential equation and the equations of motion, respectively. Even though deviations from unit length of $\vP^e_i$ do not affect the kinematic differential equation, to avoid numerical issues due to quaternion magnitudes near zero or floating point overflow, the nodal quaternions are normalized after each time-step, i.e., $\vP^e_{i} = \vP^e_i / \|\vP_i\|$. For static problems, the $n_\vu = 6N$ nonlinear generalized force equilibrium equations
\begin{equation}\label{eq:nonlinear_generalized_force_equilibrium}
\mathbf{0} =  \vf^\mathrm{int}(\vq) + \vf^\mathrm{ext}(\vq)
\end{equation}
must be augmented by the $N$ constraint equations 
\begin{equation}\label{eq:quaternion_constraints}
	\mathbf{0} = \vg(\vq) = (\|\vP^0_0\|^2 - 1, \dots, \|\vP^{n_\mathrm{el}-1}_p\|^2 - 1)
\end{equation}
to ensure solvability.

\section{Numerical experiments}\label{sec:numerical_experiments}
In the following, the quadratic strain energy density
\begin{equation}\label{eq:strain_energy_density}
W({}_K \vga, {}_K \vka_{IK}; \xi) = \frac{1}{2} \left({}_K \vga - {}_K \vga^0\right)\T \vK_{\vga} \left({}_K \vga - {}_K \vga^0\right) + \frac{1}{2} \left({}_K \vka_{IK} - {}_K \vka_{IK}^0\right)\T \vK_{\vka} \left({}_K \vka_{IK} - {}_K \vka_{IK}^0\right)
\end{equation}
is used. The superscript $0$ refers to the evaluation in the rod's reference configuration. Moreover, $\vK_{\vga} = \operatorname{diag}(EA, GA, GA)$ and $\vK_{\vka} = \operatorname{diag}(G (I_y + I_z), E I_y, E I_z)$ denote the diagonal elasticity matrices with constant coefficients given by Saint-Venant’s relations from linear elasticity. Therein, $E$ and $G$, respectively denote the Young's and shear modulus. The cross-sectional surface is denoted $A$ and $I_y$, $I_z$ are the respective second moments of area.
\subsection{Helical spring}
Following~\cite{Marino2017}, we investigate the elongation of an initially curved helical rod due to an applied external force at its tip, pointing in positive $\ve_z^I$-direction. The rod has a Young's modulus $E=10^{11}~\mathrm{N}/\mathrm{m}^2$ and Poisson's ratio $\nu=0.2$, i.e., a shear modulus $G = E / 2 (1 + \nu)$. It has an undeformed shape of a perfect helix with $n_\mathrm{c}=10$ coils, coil radius $R=10~\mathrm{mm}$, wire diameter $d=1~\mathrm{mm}$ and unloaded pitch $k=5~\mathrm{mm}$, i.e., a total height of $h=50~\mathrm{mm}$. 

In the simulation, the spring was discretized using $75$ elements of the presented finite element formulation with $p=2$. Reduced integration was performed with 2 quadrature points, while  5 points were used for all other integrals. The rod's curved initial configuration was obtained by solving the following minimization problem. Let $\xi_j = \tfrac{j}{m - 1} \in [0, 1]$ for $j \in \{0,1,\dots,m-1\}$ denote the $m$ linearly spaced evaluation points of the reference helix curve
\begin{equation}\label{eq:helix}
	\begin{aligned}
		{}_I \vr(\xi) = R \begin{pmatrix}
			\sin\varphi(\xi) \\ 
			-\cos\varphi(\xi) \\
			c \varphi(\xi)
		\end{pmatrix} 
		\, , \quad \mathrm{with} \quad c = \frac{k}{2 \pi R} \quad \mathrm{and} \quad \varphi(\xi) = 2 \pi n_\mathrm{c} \xi \, .
	\end{aligned}
\end{equation}
Hence, the evaluation of the reference curve~\eqref{eq:helix} at all $\xi_j$'s leads to $m$ target centerline points ${}_I \vr_j = {}_I \vr(\xi_j)$. Similarly, the corresponding cross-section orientations are given by evaluating the Serret--Frenet basis $\vA_{IK_j} = ({}_I \ve_x^{K_j} \ {}_I \ve_y^{K_j} \ {}_I \ve_z^{K_j})$ with ${}_I \ve_x^{K_j} = {}_I\vr_{,\xi}(\xi_j) / \|{}_I\vr_{,\xi}(\xi_j)\|$, ${}_I \ve_y^{K_j} = {}_I\vr_{,\xi\xi}(\xi_j) / \|{}_I\vr_{,\xi\xi}(\xi_j)\|$ and $\ve_z^{K_j} = {}_I \ve_x^{K_j} \times {}_I \ve_y^{K_j}$ for the individual $\xi_j$'s. Following~\cite{Harsch2023a}, the centerline positions and cross-section orientations can be assembled in the Euclidean transformations
\begin{equation}
	\vH_j = \begin{pmatrix}
		\vA_{IK_j} & {}_I \vr_{j} \\
		\mathbf{0}_{1\times3} & 1
	\end{pmatrix} \quad \mathrm{and} \quad
	\vH(\xi_j) = \begin{pmatrix}
		\vA_{IK}(\xi_j) & {}_I \vr_{OP}(\xi_j) \\
		\mathbf{0}_{1\times3} & 1
	\end{pmatrix} \, , \quad \text{with} \quad
	\vH_j^{-1} = \begin{pmatrix}
		\vA_{IK_j}^T & -\vA_{IK_j}^T \, {}_I \vr_{j} \\
		\mathbf{0}_{1\times3} & 1
	\end{pmatrix} \, .
\end{equation}
Using the $\SE(3)$-logarithm map $\operatorname{Log}_{\SE(3)}$ introduced in~\cite{Harsch2023a}, the optimal initial generalized position coordinates  $\vq_0$ results from the nonlinear least squares problem
\begin{equation}\label{eq:least_squares_SE3}
	\vq_0 = \underset{\vq \in \mathbb{R}^{n_\vq}}{\mathrm{argmin}} \, K(\vq) \, , \quad \mathrm{with} \quad K(\vq) = \frac{1}{2} \sum_{j=0}^{m-1} \|\vth_j(\vq)\|^2 \quad \mathrm{and} \quad \vth_j(\vq) = \operatorname{Log}_{\SE(3)}\big(\vH_j^{-1} \vH(\xi_j)\big) \, ,
\end{equation}
in terms of the metric of relative twists. The minimization problem~\eqref{eq:least_squares_SE3} can efficiently be solved using a Levenberg--Marquardt algorithm. The unity constraints of the nodal quaternions~\eqref{eq:quaternion_constraints} can be incorporated into the optimization problem as equality constraints, albeit at the expense of employing a complex constrained nonlinear least squares solver. To simplify the process, we initially solved the unconstrained minimization problem and subsequently applied a projection step to normalize all nodal quaternions.

\begin{figure}
	\centering
	\includegraphics{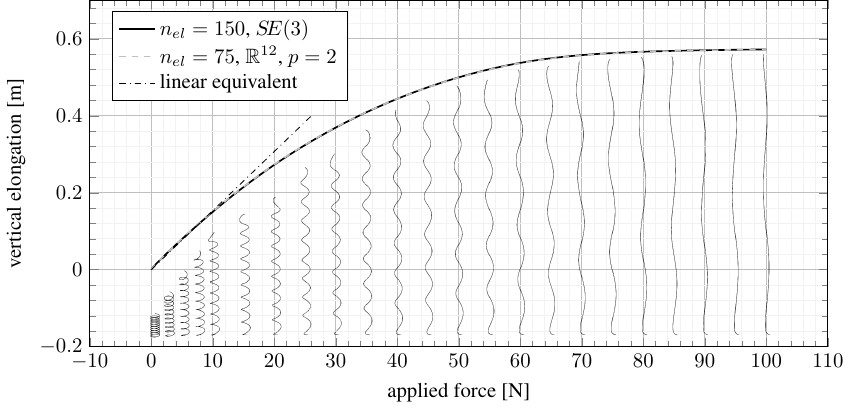}
	\caption{Force displacement diagram and deformed configurations of the helical spring.}
	\label{fig:Marino2017}
\end{figure}

Starting from $\vq_0$, the maximal force of $100~\mathrm{N}$ was applied within 500 linearly spaced force increments. During each iteration, the nonlinear equations~\eqref{eq:nonlinear_generalized_force_equilibrium} and~\eqref{eq:quaternion_constraints} were solved up to an absolute error of $10^{-8}$. As can be seen in Fig.~\ref{fig:Marino2017}, the helical spring initially elongates proportional to the applied load. This is in line with classical helical spring theory~\cite{Berg1991}, which assumes a linear force-displacement relation with linear equivalent stiffness $G d^4 / (64 n_\mathrm{c} R^3) \approx 65.1~\mathrm{N/m}$. When the elongation exceeds a certain value (approx.\,$10~\mathrm{N}$), the linear theory does not longer agree with the numerically obtained nonlinear solution. This observation was also made by~\cite{Marino2017} and can be explained as follows. The helical spring unwinds gradually and  approaches slowly a straight line with an altered linear stiffness $EA$. For comparison, we also solved the problem with the two-node $\SE(3)$-interpolation strategy proposed in~\cite{Harsch2023a}, using the same number of unknowns. As depicted in Fig.~\ref{fig:Marino2017}, the results are in line with the proposed quaternion formulation.

\subsection{Wilberforce pendulum}
\begin{figure}
	\centering
	\includegraphics{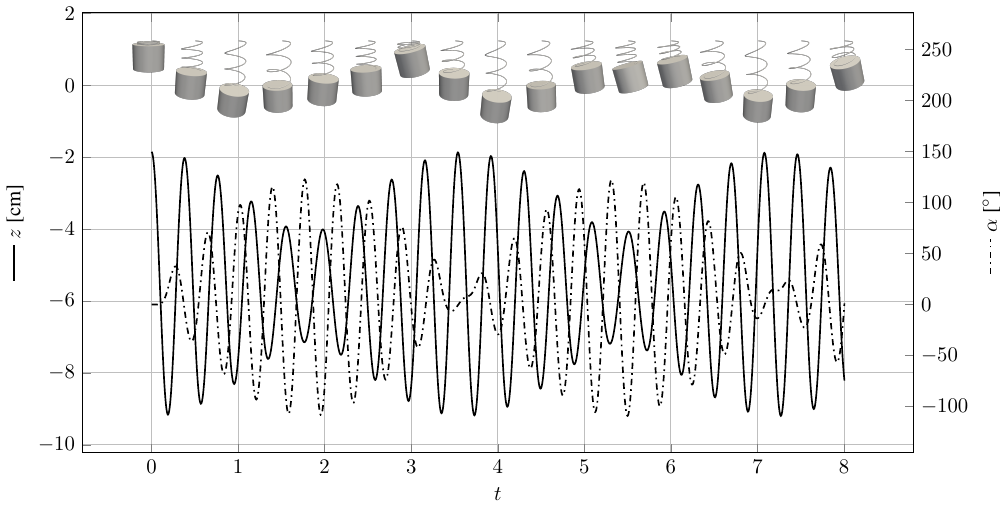}
	\caption{Vertical position $z$ of the cylinder's center of mass and rotation angle $\alpha$ corresponding to the first Euler angle with sequence ``$zyx$'' over time. Snapshots of the Wilberforce pendulum are attached to the corresponding time instants.}
	\label{fig:wilberforce}
\end{figure}

More than 100 years ago, Lionel Robert Wilberforce did investigations \emph{On the Vibrations of a Loaded Spiral Spring}~\cite{Wilberforce1894}. The experimental setup can be described as follows. While one end of a helical spring is clamped, at the other end a cylindrical bob is attached, see Fig.~\ref{fig:wilberforce}. When the cylinder in the gravitational field is displaced vertically, it starts oscillating up and down. Due to the coupling of bending and torsion of the deformed spring an additional torsional oscillation around the vertical axis of the cylinder is induced. When the cylinder's moment of inertia is properly adjusted, a beat phenomenon can be observed. In that case, the envelope of the vertical and torsional oscillations possess an almost perfect phase shift of $\pi/2$, i.e., the maximal amplitude of the vertical oscillations coincide with a zero torsional amplitude and vice versa. 

To have a benchmark example that can be reproduced with reasonably computational effort, we introduce here a Wilberforce pendulum consisting of a spring with three coils modeled as a precurved rod. The rod has the properties of steel with mass density $\rho_0=7850~\mathrm{kg}/\mathrm{m}^3$, shear modulus $G=81\cdot10^{9}~\mathrm{N}/\mathrm{m}^2$ and Poisson's ratio $\nu=0.23$, i.e., a Young's modulus $E = 2 G (1 + \nu)=199\cdot10^{9}~\mathrm{N}/\mathrm{m}^2$. The undeformed shape is given by a perfect helix with $n_\mathrm{c}=3$ coils, coil radius $R=16~\mathrm{mm}$, wire diameter $d=1~\mathrm{mm}$ and an unloaded pitch of $k=1~\mathrm{mm}$. The bob is modeled as a cylindrical rigid body with radius $r=23~\mathrm{mm}$ and height $h=36~\mathrm{mm}$ also having the mass density of steel.

In the simulations, the rod was discretized using $18$ elements of the presented Cosserat rod finite element with $p=2$. Gravitational forces for the rod were neglected. Again, reduced integration was performed with 2 quadrature points, while for all other integrals 5 points were used. The bob was parameterized by the inertial position of the center of mass ${}_I \vr_{OS}$ together with a non-unit quaternion $\vP$ for the orientation. The bob was subjected to gravity with gravity constant $g=9.81~\mathrm{m/s^2}$. For the governing equations describing such a parameterized rigid body under the influence of gravity, we refer to model 4 in~\cite{Sailer2020}. Cylinder and rod were rigidly connected by perfect bilateral constraints~\cite{Geradin2001}.
Again, the optimal helical initial configuration $\vq_0$ was found by solving the minimization problem~\eqref{eq:least_squares_SE3}. The system was initialized at rest with initial velocity $\vu_0 = \mathbf{0}$. The resulting differential algebraic equations were solved using a first-order generalized-alpha method~\cite{Jansen2000} for constrained mechanical systems of differential index 3, similar to the implementation found in~\cite{Arnold2007}. A constant step-size $\Delta t = 5\cdot10^{-3}~\mathrm{s}$ was chosen and the governing equations were solved up to a final time of $t_1 = 8~\mathrm{s}$. Since the example includes high-frequency oscillations, we chose a spectral radius at infinity of $\rho_\infty = 0.8$. The internal Newton--Raphson method satisfied a tolerance of $10^{-8}$ with respect to the maximum absolute error. In  Fig.~\ref{fig:wilberforce} the vertical position and the torsional angle of the rigid cylinder are plotted clearly showing the beat phenomenon of the Wilberforce pendulum.

\vspace{\baselineskip}


\providecommand{\WileyBibTextsc}{}
\let\textsc\WileyBibTextsc
\providecommand{\othercit}{}
\providecommand{\jr}[1]{#1}
\providecommand{\etal}{~et~al.}

\end{document}